\documentclass[12pt]{birkjour}
\usepackage{amsmath, amsthm, amscd, amsfonts,amssymb, graphicx}%
 \newtheorem{theorem}{Theorem}[section]
 \newtheorem{corollary}[theorem]{Corollary}
 \newtheorem{lemma}[theorem]{Lemma}
 \newtheorem{proposition}[theorem]{Proposition}
 \theoremstyle{definition}
 
 \newtheorem{example}[theorem]{Example}
 \theoremstyle{remark}
 
 \numberwithin{equation}{section}
\usepackage{color}
\begin{document}
\title[On Excesses and Duality in Woven Frames]
 {On Excesses and Duality in Woven Frames}

\author[Elahe Agheshteh Moghaddam]{Elahe Agheshteh Moghaddam}
\address{Department of Mathematics and Computer Sciences, Hakim Sabzevari University, Sabzevar, Iran.}
\email{elahe.moghadam@gmail.com}
\author[A. Arefijamaal]{Ali Akbar Arefijamaal}
\address{Department of Mathematics and Computer Sciences, Hakim Sabzevari University, Sabzevar, Iran.}
\email{arefijamaal@hsu.ac.ir;arefijamaal@gmail.com}


\begin{abstract} Weaving frames in separable Hilbert spaces have been recently introduced by Bemrose et al. to deal with some problems in distributed signal processing and wireless sensor networks.
In this paper, we study the notion of excess for woven frames and prove that any two frames in a separable Hilbert space
that are woven have the same excess. We also show that every frame with a large class of duals is woven provided that its redundant elements have small enough norm. Also, we try to transfer the woven property from frames to their duals and vise versa. Finally, we look at which perturbations of dual frames preserve the
woven property, moreover it is shown that under some conditions the canonical Pareseval frame of two woven frames are also woven.

 \end{abstract}

\maketitle

\section{ Introduction}
Frames which are introduced by Duffin and Schaefer provide robust, stable and usually non-unique representations of vectors \cite{Duffin}.
The theory of frames has been growing rapidly \cite{1,2,19,45} after the fundamental paper of Daubbechies, Grossmann and Meyer \cite{Daubecheis}. they have seen great achievements in pure mathematics, science, and engineering  such as  image processing, signal processing, sampling and approximation theory \cite{10,12,14,16,37,55}.

Recently, due to application and theoretical goals, some generalizations of frames have been presented, such as g-frames \cite{63}, K-frames \cite{39} and fusion frames \cite{17,mitra}. In most of them and many application problems, reconstruction and duality play a key role. Hence, characterization, construction and in general survey of duality is a valuable discussion in frame theory \cite{dual Dr,10,37}.

In signal processing each signal is interpreted as a vector. In this interpretation, a vector expressed as a linear combination of the frame elements. Using a frame, it is possible to create a simpler, more sparse representation of a signal as compared with a family of elementary signals.
The concept of $woven$ $frames$ is motivated by some problems in signal processing. For example; given two frames $\{\varphi_{i}\}_{i \in I}$ and $\{\psi_{i}\}_{i \in I}$. At each sensor we measure a signal $f$ either with $\varphi_{i}$ or with $\psi_{i}$, so that the collected information is the set of numbers $\{\langle f,\varphi_{i}\rangle\}_{i \in \sigma}\cup \{\langle f,\psi_{i}\rangle\}_{i \in \sigma^{c}}$
 for some subset $\sigma\subset I$. Now we may recover $f$ from these measurements, no matter which kind of measurement has been made at each sensors. In other words, is the set $\{\varphi_{i}\}_{i \in \sigma} \cup \{\psi_{i}\}_{i \in \sigma^{c}}$ a frame for all subset $\sigma \subset I$? This question led us to the definition of woven frames.

\section{Preliminaries and Notations}

A sequence $\Phi=\{\varphi_{i}\}_{i \in I}$ in a separable Hilbert space $\mathcal {H}$ is called a $frame$ for $\mathcal{H}$ if there exist constants $0<A_{\Phi}\leq B_{\Phi}<\infty$ such that
\begin{eqnarray}\label{frame}
A_{\Phi} \|f\|^{2}\leq \sum _{i \in I} |\langle f,\varphi_{i}\rangle|^{2} \leq B_{\Phi}\|f \|^{2}, \qquad(f \in \mathcal{H}).
\end{eqnarray}
The constants $A_{\Phi}$ and $B_{\Phi}$ are called lower and upper frame bounds, respectively. The supremum of all lower frame bounds is called the $optimal$ $lower$ $frame$ $bound$ and the infimum of all upper frame bounds is called the $optimal$ $upper$ $frame$ $bound$.
A sequence $\{\varphi_{i}\}_{i \in I}$ is called $Bessel$ if the right inequality in $(\ref{frame})$ holds.

A sequence $\Phi=\{\varphi_{i}\}_{i \in I}$ in Hilbert space $\mathcal {H}$ is called a $Riesz$ $sequence$ if there are constants $0<A_{\Phi}\leq B_{\Phi} <\infty$ so that for all finite scalars  $c_{i}$ we have
\begin{center}
$A_{\Phi} \sum_{i \in I} |c_{i}|^{2} \leq \left \| \sum _{i \in I}c_{i}\varphi_{i} \right \|^{2} \leq B_{\Phi}\sum_{i \in I} |c_{i}|^{2}.$
\end{center}
 The constants $A_{\Phi}$ and $B_{\Phi}$ are called $lower$ and $upper$ $Riesz$ $bounds$, respectively. In addition, if $\Phi$ is complete in $\mathcal{H}$, then it is a $Riesz$ $basis$ for $\mathcal{H}$.

Given a Bessel sequence $\Phi=\{\varphi_{i}\}_{i \in I}$, the $synthesis$ $operator$ $T_{\Phi} : \ell^{2}(\mathbb{N}) \rightarrow \mathcal{H}$ is defined by
$ T_{\Phi}\{c_{i}\}= \sum_{i \in I} c_{i} \varphi_{i}$. Its adjoint, $T^{*}_{\Phi}:\mathcal{H} \rightarrow \ell^{2}(\mathbb{N})$,  which is called the $analysis$ $operator$, is given by
$T^{*}_{\Phi}f =\{\langle f , \varphi_{i}\rangle\}_{i \in I}$. Moreover, $S:\mathcal{H} \rightarrow \mathcal{H}$ the $frame$ $operator$ of $\Phi$, is given by
$S_{\Phi}f=T_{\Phi}T^{*}_{\Phi}f$.
If $\Phi$ is a frame, then $S_{\Phi}$ is invertible and $A_{\Phi}\leq S_{\Phi} \leq B_{\Phi}$, see \cite{crystiansen} for more details.
The sequence $\tilde{\Phi}=\{S^{-1}_{\Phi}\varphi_{i}\}_{i \in I}$, which is also a frame, is called the $canonical$ $dual$ $frame$.
A frame $\{\psi_{i}\}_{i \in I}$ is called a $dual$ of $\{\varphi_{i}\}_{i \in I}$ if
\begin{center}
$f= \sum _{i \in I}\langle f,\psi_{i}\rangle \varphi_{i}, \qquad(f \in \mathcal{H}).$
\end{center}

\begin{theorem}\label{dual Dr} \cite{dual Dr}
If $\Phi=\{\varphi_{i}\}_{i \in I}$ be a frame. Then every dual frame of $\Phi$ is of the form of $\Phi^{d}=\{S^{-1}_{\Phi}\varphi_{i}+ u_{i}\}_{i \in I}$ where $\{u_{i}\}_{i \in I}$ is a Bessel sequence such that
\begin{eqnarray}\label{dual Drr}
 \sum_{i \in I} \langle f ,\varphi_{i}\rangle u_{i}=0, \qquad (f\in \mathcal{H}).
\end{eqnarray}
\end{theorem}

\subsection{Woven Frames}
Throughout  the paper, $\mathcal{H}$ is a separable Hilbert space, $I$ a countable index set and $I_{\mathcal{H}}$ the identity operator on Hilbert space $\mathcal{H}$. We denote $\Phi = \{\varphi_{i}\}_{i\in I}$ for a frame with bounds $A_{\Phi}$ and $B_{\Phi}$. Also we use of $[n]$ to denote the set $\{1,2, \ldots, n\}$. We now review some definitions and primary results of woven frames, which are used in the present paper. For more
information see \cite{Al00}.

A family of frames $\{\varphi_{ij}\}_{i \in I}$ for $j \in \{1,\ldots,m\}$ for a Hilbert space $\mathcal{H}$ is said to be $woven$ if there are universal constants $A$ and $B$ so that for every partition $\{\sigma_{j}\}_{j=1}^{m}$ of $I$, the family $\{\varphi_{ij}\}_{i\in \sigma_{j}, j=1}^{m}$ is a frame for $\mathcal{H}$ with lower and upper frame bounds $A$ and $B$, respectively \cite{Al00}. Each family $\{\varphi_{ij}\}_{i\in \sigma_{j}, j=1}^{m}$ is called a $weaving$.
Two frames $\{\varphi_{i}\}_{i \in I}$ and $\{\psi_{i}\}_{i \in I}$ for Hilbert space $\mathcal{H}$ are $weakly$ $woven$ if for every subset $\sigma\subset I$, the family $\{\varphi_{i}\}_{i \in \sigma}\cup \{\psi_{i}\}_{i \in \sigma^{c}}$ is a frame for $\mathcal{H}$.

Two frames $\{\varphi_{i}\}_{i \in I}$ and $\{\psi_{i}\}_{i \in I}$ are called $Riesz$ $woven$ if for every subset $\sigma\subset I$, the family $\{\varphi_{i}\}_{i \in \sigma}\cup \{\psi_{i}\}_{i \in \sigma^{c}}$ is a Riesz basis for $\mathcal{H}$. Combining Theorem 5.3 of \cite{Al00} and Theorem 3.5 of \cite{Al93}, every Riesz basis with its canonical dual is Riesz woven. However, the canonical Parseval dual of two Riesz woven bases are not necessary woven. For example, $\{e_{1}, e_{2}\}$ and $\{e_{1}+ e_{2}, 2e_{1}+e_{2}\}$ are Riesz woven bases with the canonical Parseval duals $\{e_{1}, e_{2}\}$ and $\{e_{2}, e_{1}\}$
which are clearly not woven \cite{Al00}. In this article, we study woven frames from the duality approach. One of our aim is to find dual frame pairs which are woven.

\begin{theorem}\label{lower}
If $\{\varphi_{ij}\}_{j=1, i \in I}^{M}$ is a family of Bessel sequence for $\mathcal{H}$ with a Bessel bound $B_{j}$ for all $1\leq j \leq M$ then every weaving is a Bessel sequence with the Bessel bound $\sum _{j=1}^{M}B_{j}$.
\end{theorem}

\begin{proposition}
Given two frames $\{\varphi_{i}\}_{i \in I}$ and $\{\psi_{i}\}_{i \in I}$ for $\mathcal{H}$, the following are equivalent:
\begin{itemize}
\item[(1)] The two frames are woven.
\item[(2)] The two frames are weakly woven.
\end{itemize}
\end{proposition}

\begin{proposition}\label{theo6.1}
Let $\Phi=\{\varphi_{i}\}_{i \in I}$ and $\Psi=\{\psi_{i}\}_{i \in I}$ be frames for Hilbert space $\mathcal{H}$ such that
there is a $0<\lambda<1$ so that
\begin{eqnarray*}
\lambda(\sqrt{B_{\Phi}}+\sqrt{B_{\Psi}})\leq \frac{A_{\Phi}}{2}
\end{eqnarray*}
and for all sequence of scalars $\{a_{i}\}_{i \in I}$ we have
\begin{eqnarray*}
\left\| \sum_{i \in I} a_{i}(\varphi_{i}- \psi_{i}) \right\| \leq \lambda \| \{a_{i}\}_{i \in I}\|.
\end{eqnarray*}
Then for every $\sigma \subset I$, the family $\{\varphi_{i}\}_{i \in \sigma}\cup \{\psi_{i}\}_{i \in \sigma^{c}}$ is a frame for $\mathcal{H}$  with frame bounds $\frac{A_{\Phi}}{2}$ , $B_{\Phi}+B_{\Psi}$. That is,  $\Phi$ and $\Psi$ are woven.
\end{proposition}

\begin{proposition}\label{theo6.2}
Let $\Phi=\{\varphi_{i}\}$ be a frame and $U$ be a bounded operator such that $\|I_{\mathcal{H}}-U\|^{2}<\frac{A_{\Phi}}{B_{\Phi}}$. Then $\Phi$ and $U\Phi$ are woven with the universal lower bound $\left(\sqrt{A_{\Phi}}- \sqrt{B_{\Phi}}\|I_{\mathcal{H}}-U\|\right)^{2}$.
\end{proposition}

We end this section with an example which shows that the woven property is not transitive, in general.
\begin{example}
Let $\mathcal{H}=\mathbb{R}^{2}$. Consider frames $\Phi=\{e_{1},e_{1},e_{2}\}$, $\Psi=\{e_{1},e_{2},e_{2}\}$ and $\eta=\{e_{1},e_{2},e_{1}\}$ on $\mathcal{H}$ where $\{e_{1},e_{2}\}$ is the standard orthonormal basis of $\mathcal{H}$. Then $\Phi$ with $\Psi$ and $\Psi$ with $\eta$ are woven with universal bounds $A_{1}=A_{2}=1$ and $B_{1}=B_{2}=2$. However, $\Phi$ and $\eta$ are not woven.
\end{example}

In order to achieve this, we take a restricted condition on bounds. More precisely, let $\{\varphi_{i}\}_{i \in I}$, $\{\psi_{i}\}_{i \in I}$ and $\{\eta_{i}\}_{i \in I}$ be frames for Hilbert space $\mathcal{H}$.
If $\{\varphi_{i}\}_{i \in I}$ and $\{\psi_{i}\}_{i \in I}$ are woven frames with a universal lower bound $A_{1}$, and  $\{\psi_{i}\}_{i \in I}$ is woven with $\{\eta_{i}\}_{i \in I}$  by a universal lower bound $A_{2}$ such that $A_{1}+ A_{2}- B_{\psi}>0$. Then for each $\sigma \subset I$ and $f \in \mathcal{H}$ we obtain
\begin{eqnarray*}
\left(A_{1}+ A_{2}-B_{\psi}\right)\parallel f \parallel^{2}
&\leq& \sum_{i\in \sigma}|\langle f , \varphi_{i} \rangle|^{2} + \sum_{i\in \sigma^{c}}|\langle f , \psi_{i} \rangle|^{2}
+ \sum_{i\in \sigma}|\langle f , \psi_{i} \rangle|^{2}\\
&+&\sum_{i\in \sigma^{c}}|\langle f , \eta_{i} \rangle|^{2} - \sum_ {i \in I} |\langle f , \psi_{i} \rangle|^{2} \\
&=& \sum_{i\in \sigma}|\langle f , \varphi_{i} \rangle|^{2} + \sum_{i\in \sigma^{c}}|\langle f , \eta_{i} \rangle|^{2} \\
&\leq &\left (B_{1}+ B_{2}\right)\parallel f\parallel^{2}.
\end{eqnarray*}
Hence, $\{\varphi_{i}\}_{i \in I}$ and $\{\eta_{i}\}_{i \in I}$ are woven.

\section{Weaving and Excesses}

The excess of a frame $\varphi$, denoted by $E(\varphi)$, is the greatest integer $n$ so that $n$ elements can be deleted from the frame and still leave a frame, or $+\infty$ if there is no upper bound to the number of elements that can be removed. Two frames in a separable Hilbert space $\mathcal{H}$
which are dual to each other have the same excess, see \cite{Deficit} for more details. Every frame $\varphi$ with $E(\varphi)= n$ can be written by $\varphi=\{\varphi_{i}\}_{i \in I \backslash \{i_{1}, \ldots i_{n}\}}\cup \{\varphi_{i_{1}},\ldots \varphi_{i_{n}}\}$, where $\varphi=\{\varphi_{i}\}_{i \in I \backslash \{i_{1}, \ldots i_{n}\}}$ is a Riesz basis for $\mathcal{H}$ and $\{\varphi_{i_{1}},\ldots \varphi_{i_{n}}\}$ are redundant elements of $\varphi$.

In the next theorem, we show  that woven frames have the same excess.
\begin{theorem}
Let $\Phi=\{\varphi_{i}\}_{i \in I}$ and $\Psi=\{\psi_{i}\}_{i \in I}$ be woven. Then $E(\Phi)=E(\Psi)$.
\end{theorem}
\begin{proof}
At first we suppose that $E(\Phi)=1$. Without lose of generality we can assume that
 and $\{\varphi_{i}\}_{i=2}^{\infty}$ is an orthonormal basis of $\mathcal{H}$.
If $E(\Psi)>1$ then there exist $r,s \in I$ such that $\psi_{r}, \psi_{s}\in \overline{span}\{\psi_{i}\}_{i\neq r,s}$.
Given $\epsilon > 0$ and choose $n \in \mathbb{N}$ so that
\begin{eqnarray*}
 d(\psi_{r}, K)< \epsilon, \qquad d(\psi_{s}, K)< \epsilon.
 \end{eqnarray*}
  where
\begin{eqnarray*}
 K= span \{\psi_{i}\}_{i=1, i\neq r,s}^{n}.
 \end{eqnarray*}
 Clearly, $dim K\leq n-2$ and $dim\left(span\{\varphi_{i}\}_{i=2}^{\infty}\right)=n-1$.
So, there exists $f \in span\{\varphi_{i}\}_{i=2}^{n} \cap K^{\perp}$ with $\parallel f\parallel=1$. Take $\sigma=\{1,2, \ldots, n\}$, then
\begin{eqnarray}\label{2.8}
\sum_{i \in \sigma}|\langle f, \psi_{i}\rangle |^{2} + \sum_{i \in \sigma^{c}}|\langle f, \varphi_{i}\rangle |^{2}= |\langle f, \psi_{r}\rangle |^{2} + |\langle f, \psi_{s}\rangle |^{2}.
\end{eqnarray}
For $\psi_{r}\in \mathcal{H}$ we have $\psi_{r}=\psi_{r_{1}}+ \psi_{r_{2}}$ such that $\psi_{r_{1}} \in K$ and $\psi_{r_{2}}\in K^{\perp}$. Hence,
\begin{eqnarray*}
|\langle f, \psi_{r} \rangle| &=& |\langle f, \psi_{r_{1}} \rangle+ \langle f, \psi_{r_{2}} \rangle|
= |\langle f, \psi_{r_{2}} \rangle|\\
& \leq& \|\psi_{r_{2}}\|\\
& =& \|\psi_{r}- \psi_{r_{1}}\|= d\left(\psi_{r}, K\right)< \epsilon .
\end{eqnarray*}
A similar calculation shows that $|\langle f, \psi _{s} \rangle|< \epsilon$. Therefore, by $\left(\ref{2.8}\right)$ we obtain
\begin{eqnarray*}
\sum_{i \in \sigma}|\langle f, \psi_{i}\rangle |^{2} + \sum_{i \in \sigma^{c}}|\langle f, \varphi_{i}\rangle |^{2} < 2\epsilon.
\end{eqnarray*}
So, $\{\varphi_{i}\}_{i \in \sigma} \cup \{\psi_{i}\}_{i \in \sigma^{c}}$ is not a frame and it is a contradiction. Thus, $E(\Psi)=0$ or $E(\Psi)=1$. If $E(\Psi)=0$, then $\Psi$ and $\Phi$ are Riesz bases, by Theorem $5.4$ of \cite{Al00}, which is contradiction. Therefore, $E(\Phi)=1$. By induction, the result holds for all values of $E(\Phi)=n\in \mathbb{N}$. Finally, assume that $E(\Phi)=\infty$. We can certainly see that $E(\Psi)=\infty$, since otherwise
$E(\Phi)<\infty$ which is contradiction.
\end{proof}
In \cite{Al93}, it is shown that a frame with its canonical dual is woven provided that the norm of redundant elements are small enough. In the following, we prove this fact for alternate duals.
\begin{theorem}
Let $\Phi=\{\varphi_{i}\}_{i \in I}$ be a frame for $\mathcal{H}$ such that the norm of redundant elements of $\Phi$ are small enough. Then there exist infinitely many dual frames which are woven with $\Phi$.
\end{theorem}
\begin{proof}
  Suppose first that $E(\Phi)=n$. Without lose of generality we can write $\Phi=\{\varphi_{i}\}_{i \in I\backslash [n] }\cup \{\varphi_{i}\}_{i \in [n]}$, where $\phi =\{\varphi_{i}\}_{i \in I \backslash [n]}$ is a Riesz basis for $\mathcal{H}$. Using Theorem 3.5 of \cite{Al93} follows that $\phi$ and $S^{-1}_{\phi}\phi$ are woven. In particular, $\phi$ and $S_{\phi}\phi$ are woven. Denote its universal lower bound by $A$ and suppose that $\Phi^{d}=\{S^{-1}_{\Phi}\varphi_{i}+u_{i}\}_{i\in I}$ is a dual of $\Phi$ such that $\{u_{i}\}_{i \in I}$ is a Bessel sequence satisfies in (\ref{dual Drr}). Choose $\epsilon > 0$ so that
\begin{eqnarray}\label{2.9}
\sum_{i \in [n]}\|\varphi_{i}\|^{2} +2\sqrt{\epsilon  B_{U} B_{\Phi}}<\sqrt{\frac{A}{B_{\Phi}}}.
\end{eqnarray}
 Let $\sigma\subset I$ and take $\alpha_{i,k}=\langle \varphi_{i}, \varphi_{k}\rangle, 1\leq k \leq n$. Then
\begin{eqnarray}\label{2.9.2}
\sum_{i \in \sigma}\left|\sum _{j=1}^{n}\alpha _{i,j}\langle f, \varphi_{j} \rangle\right|^{2}\leq  \parallel f \parallel ^{2}B_{\Phi}\left(\sum_{j=1}^{n}\parallel \varphi_{j} \parallel^{2}\right)^{2}, \qquad (f \in \mathcal{H}).
\end{eqnarray}
Indeed, for each $f\in \mathcal{H}$ we have
\begin{eqnarray*}
&&\sum_{i \in \sigma}\left|\sum _{j=1}^{n}\alpha _{i,j}\langle f, \varphi_{j} \rangle\right|^{2}\\
&=& \sum_{i \in \sigma} \left|\left(\alpha_{i,1}, \ldots , \alpha_{i,n}\right).\left(\overline{\langle f , \varphi_{1} \rangle }, \ldots ,\overline{\langle f , \varphi_{n} \rangle} \right)\right|^{2}\\
&\leq& \sum_{i=1}^{\infty}\left(|\alpha_{i,1}|^{2}+\ldots +| \alpha_{i,n}|^{2}\right)\left(|\langle f , \varphi_{1} \rangle|^{2}+ \ldots + |\langle f , \varphi_{n} \rangle |^{2}\right)\\
&\leq & \left(\sum _{i=1}^{\infty}|\langle \varphi_{i} , \varphi_{1} \rangle|^{2}+ \ldots + |\langle \varphi_{i} , \varphi_{n} \rangle |^{2}\right)\| f \| ^{2}\left(\| \varphi_{1} \|^{2}+\ldots + \| \varphi_{n} \| ^{2}\right)\\
&\leq& \left(B_{\Phi}\| \varphi_{1} \|^{2}+\ldots + B_{\Phi}\| \varphi_{n} \| ^{2}\right)\| f \parallel ^{2}\left(\| \varphi_{1} \parallel^{2}+\ldots + \| \varphi_{n} \| ^{2}\right)\\
&=& B_{\Phi}\left(\| \varphi_{1} \|^{2}+\ldots + \| \varphi_{n} \| ^{2}\right)^{2}\| f \| ^{2}.
\end{eqnarray*}
Putting $\Phi^{d}_{\epsilon}=\{S^{-1}_{\Phi}\varphi_{i} +\epsilon u_{i}\}_{i \in I}$. We show that $\Phi$ and $\Phi^{d}_{\epsilon}$ are woven. For this, it is enough to show the existence of a universal lower bound for $\{S_{\Phi}\varphi_{i}\}_{i\in I}$ and $\{\varphi_{i}+ S_{\Phi}\epsilon  u_{i}\}_{i\in I}$  since the woven property is preserved under bounded invertible operators \cite{Al00}. An easy computation shows that
\begin{eqnarray*}
S_{\Phi}\varphi_{i}=S_{\phi}\varphi_{i}+ \sum_{k \in [n]}\alpha_{i,k}\varphi_{k}.
\end{eqnarray*}
Let $\sigma \subset I \backslash [n]$, then we have

\begin{align*}
& \left(\sum _{i \in\sigma}|\langle f,S_{\Phi}\varphi_{i} \rangle|^{2} + \sum_{i \in \sigma^{c}}|\langle f,\varphi_{i}+ S_{\Phi}\epsilon  u_{i} \rangle|^{2}\right)^{1/2}\\
&= \Big (\sum _{i \in\sigma}\left|\langle f,S_{\Phi}\varphi_{i} \rangle\right|^{2} + \sum_{i \in\sigma^{c} \backslash [n]} \left|\langle f,\varphi_{i}+ S_{\Phi}\epsilon  u_{i} \rangle\right|^{2}+ \sum_{i \in[n]} |\langle f,\varphi_{i}+ S_{\Phi}\epsilon  u_{i} \rangle|^{2} \Big)^{1/2}\\
&=\Big (\sum _{i \in\sigma}\Big|\langle f,S_{\phi}\varphi_{i} \rangle + \sum_{k \in [n]}\alpha_{i,k}\langle f , \varphi_{k}\rangle\Big|^{2}\\
 &+ \sum_{i \in\sigma^{c} \backslash [n] }\Big|\langle f,\varphi_{i}\rangle
 + \langle f, S_{\Phi}\epsilon u_{i}\rangle\Big|^{2}
+\sum_{i \in[n]} \Big|\langle f,\varphi_{i} \rangle + \langle f, S_{\Phi}\epsilon  u_{i} \rangle\Big|^{2}\Big )^{1/2}\\
&\geq \Big(\sum _{i \in\sigma}|\langle f,S_{\phi}\varphi_{i} \rangle|^{2} + \sum_{i \in\sigma^{c} \backslash [n]} |\langle f,\varphi_{i}\rangle|^{2} + \sum_{i \in[n]} |\langle f,\varphi_{i} \rangle|^{2}\Big)^{1/2}\\
&-\left( \sum_ {i \in\sigma}\Big|\sum_{k \in [n]}\alpha_{i,k}\langle f , \varphi_{k}\rangle\Big|^{2}\right)^{1/2} - \left(\sum_{i \in\sigma^{c} \backslash [n]}| \langle f, S_{\Phi}\epsilon  u_{i} \rangle\Big|^{2}\right)^{1/2}\\ &- \left(\sum_{ [n]} |\langle f, S_{\Phi}\epsilon  u_{i} \rangle|^{2} \right)^{1/2}.
\end{align*}
Applying (\ref{2.9.2}) follows that
\begin{align*}
\sum _{i \in\sigma}|\langle f,S_{\Phi}\varphi_{i} \rangle|^{2} &+ \sum_{i \in \sigma^{c}}|\langle f,\varphi_{i}+ S_{\Phi}\epsilon  u_{i} \rangle|^{2}\\
&\geq \left(\sqrt{A}  -\sqrt{B_{\Phi}}\sum_{i \in [n]}\|\varphi_{i}\|^{2}  - 2\sqrt{\epsilon B_{U}} \parallel S_{\Phi} \parallel\right)^{2} \parallel f \parallel^{2}\\
&\geq\left( \sqrt{A} - \sqrt{B_{\Phi}}\sum_{i \in[n]}\|\varphi_{i}\|^{2}  - 2\sqrt{\epsilon B_{U}}B_{\Phi}\right)^{2}\parallel f \parallel^{2}.
\end{align*}

The desired result follows provided that $(\ref{2.9})$ holds.
As a similar calculation shows that when $\sigma^{c}\subset I \backslash [n]$ or $\sigma \subset I$, then the result is true.
 If $E(\Phi)=\infty$, we can write $\Phi =\{\varphi_{i}\}_{i \in I \backslash \sigma}\cup\{\varphi_{i}\}_{i \in \sigma}$ where $|\sigma| = \infty$ and $\phi=\{\varphi_{i}\}_{i \in I \backslash \sigma}$ is a Riesz basis for $\mathcal{H}$. Similar to the above, we can easily check that if $\sum_{i \in \sigma}\|\varphi_{i}\|^{2} 
<\sqrt{\frac{A}{B_{\Phi}}}$, where $A$ is a universal lower bound of woven frames $\Phi$ with $S_{\Phi} \Phi$, then $\Phi$ and a class of its duals are woven.
\end{proof}

 We point out that the converse of the previous theorem is not true. For instance, the frame $\Phi=\{ne_{1}, ne_{1}, e_{2}, e_{3}, \ldots\}$, where $\{e_{i}\}_{i \in I}$ is an orthonormal basis of $\mathcal{H}$ and $n\in\mathbb{N}$, is woven with all duals. Indeed, a straightforward calculation immediately shows that every dual of $\Phi$ can be written as
\begin{equation*}
\{\left(1/2n+\alpha\right)e_{1}, \left(1/2n-\alpha\right)e_{1},e_{2},e_{3}, \ldots\},
\end{equation*}
for any $\alpha\neq0$. However, the norm of redundant elements are not small enough.

\section{Weaving and Duality}
The concept of weaving not only depends on the structure of two frames, but also on the order of their elements. For example, two orthonormal bases are not necessarily woven. Thus, dual frames and the images of a frame under bounded operators are the best candidates for weaving by the original frame.
Obviously, the canonical dual of woven tight frames are woven. This leads to the following more general setting. In other word, the multiple of two woven frames are also woven.

\begin{proposition}
Let $\{\varphi_{i}\}_{i \in I}$ and $\{\psi_{i}\}_{i \in I}$ be a pair of woven frames and $\{\alpha_{i}\}_{i \in I}$ and $\{\beta_{i} \}_{i \in I}$ be a sequence of scalers such that $0 < C< |\alpha_{i}| < D < \infty$ and $0<C<  |\beta_{i}|< D< \infty$ for every $i \in I$ and for some constants C and D. Then $\{\alpha_{i} \varphi_{i}\}_{i \in I}$ and $\{\beta_{i} \psi_{i}\}_{i \in I}$ are woven.
\end{proposition}

In the next theorem, we show that if two frames $\Phi$ and $\Psi$ are woven then there are some duals $\Phi^{d}$ of $\Phi$ such that $\Psi$ and $S_{\Phi}\Phi^{d}$ are woven.


\begin{theorem}
Let $\Phi = \{\varphi_{i}\}_{i\in I}$ and $\Psi=\{\psi_i\}_{i \in I}$ be a pair of woven frames for $\mathcal{H}$ with universal bounds $A$ and $B$.
 Then there are infinitely many dual frames of $\Phi$ which are woven with $S^{-1}_{\Phi}\Psi$.
\end{theorem}
\begin{proof}
Let $U=\{u_{i}\}_{i \in I}$ be a Bessel sequence satisfies $(\ref{dual Drr})$ and take $\alpha>0$ such that
\begin{eqnarray}\label{A}
A_{\alpha}:=\alpha^{2}\|S_{\Phi}\|B_{U}+2\alpha\sqrt{B_{U}B_{\Phi}\|S_{\Phi}\|}<A.
\end{eqnarray}
 Then $\Phi_{\epsilon}^{d}=\{S^{-1}_{\Phi}\varphi_{i}+\epsilon u_{i}\}_{i \in I}$ is a dual frame of $\Phi$ for all $0<\epsilon<\alpha$.
It is enough to show that $\Psi$ and $S_{\Phi}\Phi_{\epsilon}^{d}$ are woven since the woven property is preserved under bounded invertible operators \cite{Al00}. To see this,  we only need to prove the existence of a universal lower bound.
Suppose $\sigma\subset I$, then for each $f \in \mathcal{H}$ we have

\begin{eqnarray*}
\sum_{i\in \sigma^{c}}\left|\langle f ,S_{\Phi} \Phi^{d}_{\epsilon} \rangle\right|^{2} + \sum_{i\in \sigma}\left|\langle f ,\psi_{i} \rangle \right|^{2}
&=& \sum_{i\in \sigma^{c}} \left|\langle f , \varphi_{i}+S_{\Phi}\epsilon u_{i} \rangle \right|^{2} + \sum_{i\in \sigma}|\langle f , \psi_{i} \rangle|^{2}
\\
&=&  \sum_{i\in \sigma^{c}} \left| \langle f , \varphi_{i} \rangle + \langle f , S_{\Phi}\epsilon u_{i} \rangle \right|^{2}+
 \sum_{i\in \sigma}|\langle f , \psi_{i} \rangle|^{2}\\
 & \geq & \sum_{i\in \sigma^{c}} \big||\langle f , \varphi_{i} \rangle| - | \langle f , S_{\Phi}\epsilon u_{i} \rangle| \big|^{2} +
 \sum_{i\in \sigma}|\langle f , \psi_{i} \rangle|^{2}
\\ & \geq & \sum_{i\in \sigma^{c}}\left|\langle f , \varphi_{i} \rangle\right|^{2} - \sum_{i\in \sigma^{c}} \left| \langle f , S_{\Phi}\epsilon u_{i} \rangle\right|^{2}\\
 &-&2\sum_{i\in \sigma^{c}}\left|\langle f , \varphi_{i} \rangle\right| \left| \langle f , S_{\Phi}\epsilon u_{i} \rangle\right| + \sum_{i\in \sigma}\left|\langle f , \psi_{i} \rangle\right|^{2} \\
& \geq& \left(A- \epsilon^{2}\|S_{\Phi}\|B_{U}-2\epsilon\sqrt{B_{U}}\sqrt{\|S_{\Phi}\|}\sqrt{B_{\Phi}}\right)\|f\|^{2}\\
& \geq& (A - A_{\alpha})\|f\|^{2}.
\end{eqnarray*}
By attention to $(\ref{A})$ this completes the proof.
\end{proof}

Thus, it is very natural to ask whether every frame is woven with its canonical dual. The answer is affirmative when $\mathcal{H}$ is a finite dimensional Hilbert space \cite{Al93}, however this fact is particulary true by the assumption of
 Proposition \ref{theo6.2}. Here, we obtain a sharper condition. More precisely, assume that
$\Phi = \{\varphi_{i}\}_{i \in I}$ is a frame such that
\begin{eqnarray*}
 \parallel I-S^{-1}_{\Phi}\parallel \leq \frac{A_{\Phi}}{2(B_{\Phi}+\sqrt{\frac{B_{\Phi}}{A_{\Phi}}})}.
\end{eqnarray*}
Hence,
\begin{eqnarray*}
\| T_{\Phi}-T_{\tilde{\Phi}}\|&=&\|T_{\Phi}- S^{-1}_{\Phi}T_{\Phi}\|\\
&\leq& \|T_{\Phi}\|\|I-S^{-1}_{\Phi}\|\\
&\leq&\sqrt{B_{\Phi}}\|I-S^{-1}_{\Phi}\|\\
&\leq& \frac{A_{\Phi}\sqrt{B_{\Phi}}}{2(B_{\Phi}+\sqrt{\frac{B_{\Phi}}{A_{\Phi}}})}.
\end{eqnarray*}
Combining Theorem $\ref{theo6.1}$ and Theorem $\ref{theo6.2}$ implies that  $\Phi$ and $\tilde{\Phi}$ are woven.
\begin{theorem}
Let $\Phi^{d}=\{S^{-1}_{\Phi}\varphi_{i}+ u_{i}\}_{i \in I}$ be a dual frame of $\Phi$ where $U=\{u_{i}\}_{i \in I}$ is a Bessel sequence satisfies $(\ref {dual Drr})$ such that
\begin{eqnarray}\label{2.6}
 \|I - S^{-1}_{\Phi}\| < \frac{A_{\Phi}}{2B_{\Phi}\left( 1+ \sqrt{\| S^{-1}_{\Phi} \|}\right) + 2\sqrt{B_{U}B_{\Phi}}}.
\end{eqnarray}
 Then $\Phi$  with a family of its duals are woven.
\end{theorem}
\begin{proof}
By using $(\ref{2.6})$, there exists $0<\alpha<1$ such that
\begin{eqnarray*}
\|I -S^{-1}_{\Phi}\|+\alpha \sqrt{\frac{ B_{U}}{B_{\Phi}}}
 &\leq& \frac{A_{\Phi}}{2B_{\Phi}\left(1+\sqrt{\parallel S^{-1}_{\Phi} \parallel}\right)+2\sqrt{B_{U}B_{\Phi}}}\\
&\leq& \frac{A_{\Phi}}{2B_{\Phi}\left(1+\sqrt{\parallel S^{-1}_{\Phi} \parallel}\right)+2\alpha\sqrt{ B_{U}B_{\Phi}}}.\\
\end{eqnarray*}
Putting $\Phi_{\alpha}^{d}=\{S^{-1}_{\Phi}\varphi_{i}+\alpha u_{i}\}_{i \in I}$. Then $\Phi_{\alpha}^{d}$ is a dual frame of $\Phi$ and an easy calculation follows that
\begin{eqnarray*}
B_{\Phi_{\alpha}^{d}}=\left(\sqrt{\parallel S^{-1}_{\Phi}\parallel B_{\Phi}}+\alpha\sqrt{B_{U}}\right)^{2}.
\end{eqnarray*}
 Using this fact and the above inequality we obtain
\begin{eqnarray*}
\|T_{\Phi} -T_{\Phi^{d}_{\alpha}}\| &=&\|T_{\Phi}- S^{-1}_{\Phi}T_{\Phi} -\alpha T_{U}\|\\
&\leq&\|I - S^{-1}_{\Phi}\|\|T_{\Phi}\| +\alpha \|T_{U}\|\\
&\leq&\|I - S^{-1}_{\Phi}\|\sqrt{B_{\Phi}} + \alpha  \sqrt{B_{U}}\\
&\leq& \sqrt{B_{\Phi}}  \frac{A_{\Phi}}{2B_{\Phi}\left(1+\sqrt{\parallel S^{-1}_{\Phi} \parallel}\right)+2\alpha\sqrt{B_{U}B_{\Phi}}}\\
&=&  \frac{A_{\Phi}}{2\left(\sqrt{B_{\Phi}}+\sqrt{B_{\Phi}} \sqrt{\parallel S^{-1}_{\Phi} \parallel}+\alpha\sqrt{B_{U}}\right)}\\
&=& \frac{A_{\Phi}}{2\left(\sqrt{B_{\Phi}} + \sqrt{B_{\Phi_{\alpha}^{d}}} \right)}.
\end{eqnarray*}
Applying Theorem $\ref{theo6.1}$ follows that $\Phi$ and $\Phi_{\alpha}^{d}$ are woven.
\end{proof}

In the next theorem, we show that the wovenness property can be transferred to special dual frames. First, we prove it for Parseval frames.
\begin{lemma}
Let $\Phi=\{\varphi_{i}\}_{i \in I}$ and $\Psi=\{\psi_{i}\}_{i \in I}$ be Parseval woven frames. Then there are infinitely many dual frames $\Phi^{d}$ of $\Phi$ and $\Psi^{d}$ of $\Psi$ such that $\Phi^{d}$ and $\Psi^{d}$ are woven.
\end{lemma}
\begin{proof}
Using Theorem \ref{dual Dr} follows that dual frames of $\Phi$ and $\Psi$ are of the form $\{\varphi_{i} + u_{i}\}_{i \in I}$ and $\{\psi_{i} + v_{i}\}_{i \in I}$, respectively, where  $\{u_{i}\}_{i \in I}$ and $\{v_{i}\}_{i \in I}$  are Bessel sequences and for every $f\in \mathcal{H}$
\begin{eqnarray*}
 \sum_{i \in I} \langle f ,\varphi_{i}\rangle u_{i}=0, \qquad \sum_{i \in I} \langle f ,\psi_{i}\rangle v_{i}=0.
\end{eqnarray*}
Choose $\alpha_{0} >0$ such that
\begin{eqnarray}
\alpha_{0}  B_{U} +\alpha_{0}  B_{V} + 2 \sqrt{\alpha_{0}  B_{U} B_{\Phi}} + 2\sqrt{\alpha_{0} B_{V} B_{\Psi}}<A,
\end{eqnarray}
where $A$ is a universal lower bound of weaving $\Phi$ and $\Psi$. Hence, for every $\sigma \subset I$ we have
\begin{eqnarray*}
   \sum_{\sigma}|\langle f,\varphi_{i} +\alpha_{0} u_{i} \rangle|^{2}&+&\sum_{\sigma^{c}}|\langle f,\psi_{i} +\alpha_{0} v_{i}\rangle |^{2}\\
 &=& \sum_{\sigma}|\langle  f,\varphi_{i}\rangle +  \langle f,\alpha_{0}  u_{i} \rangle |^{2}+\sum_{\sigma^{c}}|\langle  f,\psi_{i} \rangle + \langle f, \alpha_{0}  v_{i}\rangle |^{2} \\
 &\geq& \sum_{\sigma}\Big| |\langle  f,\varphi_{i}\rangle| - |\langle f,\alpha_{0}  u_{i} \rangle| \Big|^{2}+\sum_{\sigma^{c}}\Big||\langle  f,\psi_{i} \rangle| - |\langle f,\alpha_{0}  v_{i}\rangle| \Big|^{2} \\
  &\geq& \sum_{\sigma}|\langle  f,\varphi_{i}\rangle|^{2} - \sum_{\sigma}|\langle f,\alpha_{0}  u_{i} \rangle|^{2} -2  \sum_{\sigma}|\langle  f,\varphi_{i}\rangle||\langle f,\alpha_{0}  u_{i} \rangle|\\
 &+& \sum_{\sigma^{c}}|\langle  f,\psi_{i} \rangle|^{2} -  \sum_{\sigma^{c}}|\langle f,\alpha_{0}  v_{i}\rangle|^{2 } -2  \sum_{\sigma^{c}}|\langle  f,\psi_{i}\rangle||\langle f,\alpha_{0}  v_{i} \rangle|\\
 &\geq& \left(A -\alpha_{0} B_{u} -\alpha_{0} B_{v} -2 \sqrt{\alpha_{0}  B_{u} B_{\varphi}} -2\sqrt{\alpha_{0}  B_{v} B_{\psi}} \right)\|f\|^{2}
\end{eqnarray*}
This shows that $\Phi^{d}_{\alpha}=\{\varphi_{i} +\alpha u_{i}\}_{i \in I}$ and $\Psi^{d}_{\alpha}=\{\psi_{i} +\alpha v_{i}\}_{i \in I}$ which are dual frames of $\Phi$ and $\Psi$, respectively, are woven for all $0< \alpha < \alpha_{0}$.
\end{proof}

\section{Stability of dual woven frames}
In this section, we state some general stability results for woven frames, compare to the facts in \cite{Al93,Bala19}.

\begin{theorem}
Let $\Phi=\{\varphi_{i}\}_{i \in I}$ and $\Psi=\{\psi_{i}\}_{i \in I}$ be a pair of woven frames such that for every finite scalar sequence $\{c_i\}$
\begin{eqnarray}\label{AB}
\left\|\sum_{i \in I}c_{i}\varphi_{i}- c_{i}\psi_{i}\right\| < \frac{\sqrt{A}}{\sqrt{B_{\Psi}}\|S^{-1}_{\Phi}\|\|S^{-1}_{\Psi}\|\left(\sqrt{B_{\Psi}}+\sqrt{B_{\Phi}}\right)}\sum_{i \in I}|c_{i}|^{2},
\end{eqnarray}
where $A$ is a universal lower bound of weaving $S^{-1}_{\Phi}\Phi$ and $S^{-1}_{\Phi}\Psi$. Then there are infinitely many dual frames $\Phi^{d}$ of $\Phi$ and $\Psi^{d}$ of $\Psi$ such that $\Phi^{d}$ and $\Psi^{d}$ are woven.
\end{theorem}
\begin{proof}
By the assumption $S^{-1}_{\Phi}\Phi$ and $S^{-1}_{\Phi}\Psi$  are also woven. Denote its universal  bounds by $A$ and $B$.
Choose arbitrary dual frames $\{S^{-1}_{\Phi}\varphi_{i} +  u_{i}\}_{i \in I}$ and $\{S^{-1}_{\Psi}\psi_{i} + v_{i}\}_{i \in I}$ where  $U=\{u_{i}\}_{i \in I}$ and $V=\{v_{i}\}_{i \in I}$ are Bessel sequences satisfying (\ref{dual Drr}), see Theorem \ref{dual Dr}. Then
\begin{eqnarray*}
\|S_{\Phi}- S_{\Psi}\|&=& \|T_{\Phi}T^{*}_{\Phi} - T_{\Psi}T^{*}_{\Psi}\|= \|T_{\Phi}T^{*}_{\Phi} -T_{\Phi}T^{*}_{\Psi} + T_{\Phi}T^{*}_{\Psi} - T_{\Psi}T^{*}_{\Psi}\|\\
&\leq& \|T_{\Phi}- T_{\Psi}\|\left(\|T_{\Phi}\|+\|T_{\Psi}\|\right).
\end{eqnarray*}
 By (\ref{AB}), we obtain
\begin{eqnarray*}
\|S_{\Phi}- S_{\Psi}\|<\frac{1}{\left\|S^{-1}_{\Phi}\right\|\left\|S^{-1}_{\Psi}\right\|} \sqrt{\frac{A}{B_{\Psi}}}.
\end{eqnarray*}
So,
\begin{eqnarray}\label{abb}
\|S^{-1}_{\Phi} - S^{-1}_{\Psi} \| = \| S^{-1}_{\Psi} \left(S_{\Phi}-S_{\Psi} \right)S^{-1}_{\Phi}\|
< \sqrt{\frac{A}{B_{\Psi}}}.
\end{eqnarray}
Choose $0< \alpha <1$ such that
\begin{eqnarray}\label{123}
\|S_{\Phi}^{-1}-S_{\Psi}^{-1}\|+\alpha \left(\frac{\sqrt{B_{U}}+ \sqrt{B_{V}}}{\sqrt{B_{\Psi}}}\right)< \sqrt{\frac{A}{B_{\Psi}}}.
\end{eqnarray}
Then $\Phi^{d}_{\alpha}=\{S^{-1}_{\Phi}\varphi_{i} + \alpha u_{i}\}_{i \in I}$ and $\Psi^{d}_{\alpha}=\{S^{-1}_{\Psi}\psi_{i} +\alpha v_{i}\}_{i \in I}$ are duals of $\Phi$ and $\Psi$, respectively . Hence, for every $\sigma \subset I$ we have
\begin{align*}
& \Big(\sum_{\sigma}|\langle f,\Phi^{d}_{\alpha} \rangle |^{2}+\sum_{\sigma^{c}}|\langle  f,\Psi^{d}_{\alpha}\rangle |^{2}\Big)^{\frac{1}{2}}\\
&= \Big(\sum_{\sigma}\Big|\langle S^{-1}_{\Phi} f,\varphi_{i}\rangle + \langle f,\alpha  u_{i} \rangle \Big|^{2}+\sum_{\sigma^{c}}\Big|\langle S^{-1}_{\Psi} f,\psi_{i} \rangle + \langle f,\alpha  v_{i}\rangle \Big|^{2}\Big)^{\frac{1}{2}} \\
&= \Big(\sum_{\sigma}\Big|\langle S^{-1}_{\Phi} f,\varphi_{i}\rangle + \langle f,\alpha  u_{i} \rangle \Big|^{2} + \sum_{\sigma^{c}}\Big|\langle S^{-1}_{\Phi} +(S^{-1}_{\Psi}-S^{-1}_{\Phi}) f,\psi_{i} \rangle + \langle f, \alpha  v_{i}\rangle \Big|^{2}\Big)^{\frac{1}{2}} \\
&= \Big(\sum_{\sigma}\Big|\langle S^{-1}_{\Phi} f,\varphi_{i}\rangle + \langle f, \alpha  u_{i} \rangle \Big|^{2} + \sum_{\sigma^{c}}\Big|\langle S^{-1}_{\Phi}f, \psi_{i}\rangle + \langle(S^{-1}_{\Psi}-S^{-1}_{\Phi}) f,\psi_{i} \rangle + \langle f, \alpha  v_{i}\rangle \Big|^{2}\Big)^{\frac{1}{2}} \\
 &\geq \Big(\sum_{\sigma}|\langle S^{-1}_{\Phi} f,\varphi_{i}\rangle|^{2} + \sum_{\sigma^{c}}|\langle S^{-1}_{\Phi}f, \psi_{i}\rangle|^{2}\Big)^{\frac{1}{2}} -\Big(\sum_{\sigma}|\langle f,\alpha  u_{i} \rangle |^{2}\Big)^\frac{1}{2} - \Big(\sum_{\sigma^{c}}|\langle f,\alpha  v_{i} \rangle |^{2}\Big)^{\frac{1}{2}}\\
 &\quad - \Big(\sum_{\sigma^{c}}|\langle(S^{-1}_{\Psi}-S^{-1}_{\Phi}) f,\psi_{i} \rangle|^{2}\Big)^{\frac{1}{2}}\\
&\geq\left( \sqrt{A} - \alpha \sqrt{B_{U}} - \alpha \sqrt{B_{V}}- \sqrt{B_{\Psi}}\parallel S^{-1}_{\Psi}-S^{-1}_{\Phi}\parallel\right) \parallel f\parallel\\
&\geq \sqrt{B_{\Psi}}\left(\sqrt{\frac{A}{B_{\Psi}}} - \alpha \left(\frac{\sqrt{B_{U}} + \sqrt{B_{V}}}{\sqrt{B_{\Psi}}} \right) - \parallel S^{-1}_{\Psi}-S^{-1}_{\Phi}\parallel \right)\parallel f\parallel.
\end{align*}
Applying $(\ref{123})$, we obtain a universal lower bound for $\Phi^{d}_{\alpha}$ and $\Psi^{d}_{\alpha}$. This completes the proof.
\end{proof}
The next result, which is proved similarly, shows that two frames are woven from the wovenness of their duals.
\begin{corollary}
Let $\Phi=\{\varphi_{i}\}_{i \in I}$ and $\Psi=\{\psi_{i}\}_{i \in I}$ be frames with dual frames $\Phi^{d}= \{S^{-1}_{\Phi}\varphi_{i} + u_{i}\}_{i \in I}$ and $\Psi^{d}=\{S^{-1}_{\Psi}\psi_{i} + v_{i}\}_{i \in I}$, respectively. If $\Phi^{d}$ and $\Psi^{d}$ are woven with universal bounds $A$ and $B$ such that
\begin{eqnarray}\label{2.n}
 \sqrt{B_{U}}  +  \sqrt{B_{\Psi}}\parallel S_{\Phi^{d}}-S_{\Psi^{d}}\parallel< \sqrt{A}.
\end{eqnarray}
Then $\Phi$ and $\Psi$ are woven.
\end{corollary}
Now, let us restrict our attention to the canonical duals, we obtain the same conditions.
\begin{corollary} Let $\Phi$ and $\Psi$ be two frames. The following assertions are hold:
\begin{itemize}
\item[(1)] If $\Phi$ and $\Psi$ are woven frames with a universal lower bound $A$ such that
\begin{align*}
\parallel S ^{-1}_{\Phi}-  S ^{-1}_{\Psi}\parallel \leq \sqrt{\frac{A}{B_{\Psi}}}.
\end{align*}
Then $\tilde{\Phi}$ and $\tilde{\Psi}$ are woven frames.
\item[(2)] If  $\tilde{\Phi}$ and $\tilde{\Psi}$ are woven frames with a universal lower bound $A$ such that
\begin{align*}
\parallel S ^{-1}_{\Phi}-  S ^{-1}_{\Psi}\parallel \leq \sqrt{\frac{A}{B_{\Psi}}}.
\end{align*}
Then $\Phi$ and $\Psi$ are woven frames.
\end{itemize}
\end{corollary}

For every frame $\Phi=\{\varphi_{i}\}_{i \in I}$, the Parseval frame $\{S_{\Phi}^{\frac{-1}{2}}\varphi_{i}\}_{i \in I}$ is called the canonical Parseval frame of $\Phi$. In \cite{Al00}, it is represented an example that two frames are woven but their canonical Parseval frames are not woven. In the next theorem, we show that this fact is true under some condition.

\begin{theorem}
Let $\Phi$ and $\Psi$ be a pair of woven frames with a universal lower bound $A$ such that
\begin{eqnarray}\label{canonical parseval}
\|S_{\Phi}^{\frac{1}{2}} - S_{\Psi}^{\frac{1}{2}}\|< \frac{1}{\|S_{\Phi}^{-\frac{1}{2}}\|\|S_{\Psi}^{-\frac{1}{2}}\|}\left( \frac{ \sqrt{\|S^{-\frac{1}{2}}_{\Psi}\|^{2}+ A B_{\Phi}}-\|S^{-\frac{1}{2}}_{\Psi}\|}{B_{\Phi}}\right).
\end{eqnarray}
Then their canonical Parseval frames $\Phi$ and $\Psi$ are also woven.
\end{theorem}
\begin{proof} It is not difficult to see that
\begin{eqnarray}\label{canonical parseval3}
B_{\Phi} \left\|S^{-\frac{1}{2}}_{\Phi} - S^{-\frac{1}{2}}_{\Psi}\right \|^{2}  +2 \left\|S^{-\frac{1}{2}}_{\Phi} - S^{-\frac{1}{2}}_{\Psi}\right \| \left\| S^{-\frac{1}{2}}_{\Psi} \right\| - A<0,
\end{eqnarray}
whenever we have
\begin{eqnarray*}
\left\|S^{-\frac{1}{2}}_{\Phi} - S^{-\frac{1}{2}}_{\Psi}\right \|
< \frac{-\|S^{-\frac{1}{2}}_{\Psi}\|+ \sqrt{\|S^{-\frac{1}{2}}_{\Psi}\|^{2}+ A B_{\Phi}}}{B_{\Phi}}.
\end{eqnarray*}
Also, by the assumption we have
\begin{eqnarray*}
 \left\|S^{-\frac{1}{2}}_{\Phi} - S^{-\frac{1}{2}}_{\Psi}\right \|
  &=&  \left\|S^{-\frac{1}{2}}_{\Phi}\left( S^{\frac{1}{2}}_{\Psi} - S^{\frac{1}{2}}_{\Phi}  \right)  S^{-\frac{1}{2}}_{\Psi}\right \|\\
 &\leq& \left\|S^{-\frac{1}{2}}_{\Phi} \right \| \left\| S^{\frac{1}{2}}_{\Psi} - S^{\frac{1}{2}}_{\Phi}\right\| \left\| S^{-\frac{1}{2}}_{\Psi}\right \|\\
(\textrm{Due to (\ref{canonical parseval})})&<& \left\|S^{-\frac{1}{2}}_{\Phi} \right \|\left\| S^{-\frac{1}{2}}_{\Psi}\right \|\frac{1}{\|S_{\Phi}^{-\frac{1}{2}}\|\|S_{\Psi}^{-\frac{1}{2}}\|}\left( \frac{\sqrt{\|S^{-\frac{1}{2}}_{\Psi}\|^{2}+ A B_{\Phi}}-\|S^{-\frac{1}{2}}_{\Psi}\|}{B_{\Phi}}\right)\\
&=& \frac{-\|S^{-\frac{1}{2}}_{\Psi}\|+ \sqrt{\|S^{-\frac{1}{2}}_{\Psi}\|^{2}+ A B_{\Phi}}}{B_{\Phi}}.
\end{eqnarray*}
So,
\begin{eqnarray*}
\left\|S^{-\frac{1}{2}}_{\Phi} - S^{-\frac{1}{2}}_{\Psi}\right \|  < \frac{-\|S^{-\frac{1}{2}}_{\Psi}\|+ \sqrt{\|S^{-\frac{1}{2}}_{\Psi}\|^{2}+ A B_{\Phi}}}{B_{\Phi}}.
\end{eqnarray*}
In particular, $(\ref{canonical parseval3})$ is valid.
Since $\Phi$ and $\Psi$ are woven so for every $\sigma \subset I$ we have
\begin{align*}
\sum_{\sigma}|\langle  f,S^{-\frac{1}{2}}_{\Phi} \varphi_{i} \rangle |^{2}&+\sum_{\sigma^{c}}|\langle f,S^{-\frac{1}{2}}_{\Psi}\psi_{i}\rangle |^{2}\\
 &= \sum_{\sigma}\left|\langle S^{-\frac{1}{2}}_{\Phi} f, \varphi_{i} \rangle \right|^{2}+\sum_{\sigma^{c}}\left|\langle S^{-\frac{1}{2}}_{\Psi} f,\psi_{i}\rangle \right|^{2}\\
 &= \sum_{\sigma}\left|\langle\left( S^{-\frac{1}{2}}_{\Phi} - S^{-\frac{1}{2}}_{\Psi} \right) f, \varphi_{i} \rangle +  \langle S^{-\frac{1}{2}}_{\Psi} f,\varphi_{i} \rangle \right|^{2}+\sum_{\sigma^{c}}|\langle S^{-\frac{1}{2}}_{\Psi} f,\psi_{i}\rangle |^{2}\\
 &\geq \sum_{\sigma}\left|\langle S^{-\frac{1}{2}}_{\Psi} f, \varphi_{i} \rangle \right|^{2}+\sum_{\sigma^{c}}\left|\langle S^{-\frac{1}{2}}_{\Psi}f , \psi_{i} \rangle\right|^{2} - \sum_{\sigma}\left|\langle \left( S^{-\frac{1}{2}}_{\Phi} - S^{-\frac{1}{2}}_{\Psi} \right) f, \varphi_{i} \rangle \right|^{2}\\
 &\qquad -2 \sum_{\sigma}\left|\langle S^{-\frac{1}{2}}_{\Psi} f, \varphi_{i} \rangle \right|\left|\langle \left( S^{-\frac{1}{2}}_{\Phi} - S^{-\frac{1}{2}}_{\Psi} \right) f, \varphi_{i} \rangle \right|\\
 &\geq\left( A - B_{\Phi} \left\|S^{-\frac{1}{2}}_{\Phi} - S^{-\frac{1}{2}}_{\Psi}\right \|^{2}  -2 \left\|S^{-\frac{1}{2}}_{\Phi} - S^{-\frac{1}{2}}_{\Psi}\right \| \left\| S^{-\frac{1}{2}}_{\Psi} \right\| \right)\|f\|^{2}.
\end{align*}
By attention to $(\ref{canonical parseval3})$, Parseval frames $\{S^{-\frac{1}{2}}_{\Phi}\varphi_{i}\}_{i \in I}$ and $\{S^{-\frac{1}{2}}_{\Psi}\psi_{i}\}_{i \in I}$ are woven.
\end{proof}


\begin{thebibliography}{99}

\bibitem{1}
A.~Aldroubi,
\newblock Portraits of Frames,
\newblock{\em Proc. Amer. Math. Soc,} \textbf{123} (1995), 1661--1668.

\bibitem{2}
S.~T.~Ali, J.~P.~Antoine and J.~P.~ Gazean,
\newblock Continous Frames in Hilbert Spaces,
\newblock{\em Ann. Physic,} \textbf{222} (1993), 1--37.

\bibitem{Al93}
F.~Arabyani Neyshaburi and A.~Arefijamaal,
\newblock Manufactoring Pairs of Woven Frames Applying Duality Principle on Hilbert Spaces,
\newblock {\em Bulletin of the Malaysian,}  (2020).

\bibitem{10}
A.~Arefijamaal and E.~Zekaee,
\newblock Image Processing by Alternate Dual Gabor Frames,
\newblock {\em Bull. Iranian Math. Soc,} \textbf{42}(6) (2016), 1305--1314.

\bibitem{dual Dr}
A.~Arefijamaal and E.~Zekaee,
\newblock Signal Processing by Alternate Dual Gabor Frames,
\newblock{\em Appl. Comput. Harmon. Anal,} \textbf{35} (2013), 535--540.

\bibitem{Deficit}
D.~Baki\'{c}, t.~Beri\'{c},
 \newblock On excesses of frames,
 \newblock{\em Glas. Mat. Ser. III}, \textbf{50}(70), 415-427.


\bibitem{Bala19}
P.~Balazs, M.~Shamsabadi, A.~Arefjamaal, and A.~Rahimi,
\newblock U-cross Gram matrices and their invertibility.
\newblock {\em J. Math. Anal. Appl,} \textbf{476}(2) (2019), 367-390.

\bibitem{Al00}
 T.~Bemrose, P. G.~Casazza, K.~Grochenig,M. C. Lammers and R. G ~Lynch,
\newblock Weaving Hilbert space frames,
\newblock{\em Operators and Matrices,} \textbf{10}(4) (2016), 1093-1116.

\bibitem{12}
J.~Benedeto, A.~Powell and O.~Yilmaz,
\newblock Sigm-Delta Quantization and Finite Frames,
\newblock {\em IEEE Trans. Infrm. Theory,} \textbf{52} (2006), 1990-2005.

\bibitem{14}
H.~Boleskel,F.~Hlawatsch and H. G.~Feichtinger ,
\newblock Frae Theoretic Analysis of Oversampled Filter Banks,
\newblock{\em  IEEE Trans. Signal Process,} \textbf{46} (1998), 3256-3268.

\bibitem{17}
P. G.~Cazass and  G.~Kutyniok,
\newblock Frames of Subspaces,
\newblock{\em Contempt . Math}, \textbf{345} (2004), 87-114.


\bibitem{19}
 P. G.~Casazza,
\newblock The art of Frame Theory,
\newblock {\em Taiwanese Jour. Math,} \textbf{4}(2) (2000), 129--202.


\bibitem{16}
P. G.~Casazza, G.~Kutyniok, S.~Li and C.J.~ Rozell,
\newblock Modeling Sensor Networks with Fusion Frames,
\newblock {\em Wavelets XII. San Diego , CA,}
\newblock{\em  SPIE Proc. SPIE, Bellingham, WA}, (2007), 67011M-1-67011M-11.


\bibitem{crystiansen}
O.~Christensen,
\newblock Frames and Bases: An Introductory Course,
\newblock {\em Birkh\"{a}user, Boston}, (2008).

\bibitem{Daubecheis}
I.~Daubecheis, A.~Grossmann and Y.~ Meyer,
\newblock  Painless Nonorthogonal Expansions,
\newblock {\em J. Math. Phys}, \textbf{27} (1986), 1271-1283.

\bibitem{Duffin}
R.~Duffin and A.~Schaeffer,
\newblock A Class of Nonharmonic Fourier Series,
\newblock {\em Trans. Amer. Math. soc}, \textbf{72} (1952), 341-366.



\bibitem{37}
M.~Gaianu and D.M.~Onchis,
\newblock Face and Marker Detection Using Gabor Frames on GPUs,
\newblock {\em  Signal Processing}, \textbf{96} (2014), 90- 93.

\bibitem{39}
L.~G\v{a}vru\c{t}a,
\newblock Frames for Operators,
\newblock{\em Appl. Comp. Harm. Anal. Appl}, \textbf{32} (2012), 139-144.

\bibitem{45}
A. J. E. M.~Janssen,
\newblock Duality and Biorthogonality for Weyl-Heisenberg Frames,
\newblock {\em J. Fourier Anal. Appl}, \textbf{1}(4) (1995), 403-436.

\bibitem{55}
S.~Mallat,
\newblock A Wavelet Tour of Signal Processing,
\newblock {\em Academic Press, second edition}, (1999).



\bibitem{mitra}
M.~Shamsabadi and A.~A. Arefijamaal.
\newblock The invertibility of fusion frame multipliers,
\newblock {\em Linear and Multilinear Algebra}, \textbf{65}(5) (2016), 1062--1072.


\bibitem{63}
W.~Sun,
\newblock G-frames and G- Riesz bases,
\newblock {\em J.Math. Anal. Appl}, \textbf{322} (2006), 437-452.

\end{thebibliography}
\end{document}